\documentclass[12pt]{article}%
\usepackage[applemac]{inputenc}
\usepackage{amsmath,amssymb}
\usepackage{amsfonts}
\usepackage{amsmath,amssymb}
\usepackage[applemac]{inputenc}
\usepackage{color}
\usepackage{color, colortbl}
\usepackage{amsmath}
\usepackage{amssymb}
\usepackage{graphicx}
\usepackage{tikz}
\usepackage{geometry}
\usetikzlibrary{arrows}
\usepackage{amsfonts}
\usepackage{amsmath,amssymb}
\usepackage[applemac]{inputenc}
\usepackage{color}
\usepackage{amsmath}
\usepackage{amssymb}
\usepackage{graphicx}
\usepackage{tikz}
\newtheorem{theorem}{Theorem}[section]
\newtheorem{lemma}[theorem]{Lemma}

\newtheorem{definition}[theorem]{Definition}

\newtheorem{example}[theorem]{Example}

\newtheorem{problem}[theorem]{Problem}
\newtheorem{remark}[theorem]{Remark}

\usetikzlibrary{arrows,shapes,automata,petri}

\newcommand{\dproof}{\noindent {Proof.} \quad}
\newcommand{\fproof}{\hfill $\square$ \bigskip}

\numberwithin{equation}{section}

\definecolor{LightCyan}{rgb}{0.88,1,1}

\def\1B{\text{1\!\!I}}

\def\<{\langle}
\def\>{\rangle}

\def\tilde{\widetilde}
\def\hat{\widehat}
\def\R{\mathbb{R}}
\begin{document}

\title{Optimal control of SPDEs driven by time-space Brownian motion}
\author{Nacira Agram$^{1},$ Bernt \O ksendal$^{2}$, Frank Proske$^{2}$ and Olena Tymoshenko$^{2,3}$}
\date{Corrected version 31 December 2023}
\maketitle

\footnotetext[1]{Department of Mathematics, KTH Royal Institute of Technology 100 44, Stockholm, Sweden. \newline
Email: nacira@kth.se. Work supported by the Swedish Research Council grant (2020-04697).}

\footnotetext[2]{%
Department of Mathematics, University of Oslo, Norway. \\
Emails: oksendal@math.uio.no, proske@math.uio.no, otymoshenkokpi@gmail.com}

\footnotetext[1]{Department of Mathematics, KTH Royal Institute of Technology 100 44, Stockholm, Sweden. \newline
Email: nacira@kth.se. Work supported by the Swedish Research Council grant (2020-04697).}

\footnotetext[2]{%
Department of Mathematics, University of Oslo, Norway. \\
Emails: oksendal@math.uio.no, proske@math.uio.no, otymoshenkokpi@gmail.com}
\footnotetext[3]
{Department of Mathematical Analysis and Probability Theory, NTUU Igor Sikorsky Kyiv Polytechnic Institute, Kyiv, Ukraine.}

\vskip 0.6cm
\textcolor{blue}{Dedicated to the memory of Yuri Kondratiev}

\begin{abstract}
In this paper we study a Pontryagin type stochastic maximum principle for the optimal control of a system, where the state dynamics satisfy a stochastic partial differential equation (SPDE)
driven by a two-parameter (time-space) Brownian motion (also called Brownian sheet). 

We first discuss some properties of a Brownian sheet driven linear SPDE which models the growth of an ecosystem. 

Further, applying time-space white noise calculus we derive sufficient conditions and necessary conditions of optimality of the control. 

Finally, we illustrate our results by solving a linear quadratic control problem and an optimal harvesting problem in the plane. We also study possible applications to machine learning.
\end{abstract}

%\tableofcontents
%BEGINNING
%\selectlanguage{english}
%%%%%%%%%%%%%%%%%%%%%%%%%%%%%%%%%%%%%%%%%%%%%%%%%%%%%%%%%%%%%%%%%%%%%%%%%%%%%%%%%%%%%%%%%%%%%%%%%%%%%%%%%%%%%%%%%%%%%%%%

%\textbf{Keywords:} SPDE, two-parameter Brownian motion, optimal control, maximum principles, BSPDE in the plane, linear-quadratic control, machine learning.  
%%%%%%%%%%%%%%%%%%%%%%%%%%%%%%%%%

\section{Introduction}
%\frametitle{1.Introduction}
The purpose of this paper is to study optimal control of systems driven by the Brownian sheet. \vskip 0.2cm
Throughout this work, we denote by $\{B(t,x): t \geq 0, x \in \mathbb{R}\}$ a  Brownian sheet and $(\Omega, \mathcal{F}, P)$ a complete probability space on which we define the (completed) $\sigma$-field $\mathcal{F}_{t,x}$ generated by $ B(s,a), s \leq t, a \leq x$.

Wong \& Zakai \cite{WZ} generalised the notion of stochastic integrals with respect to 1-parameter Brownian motion to stochastic integrals driven by the two-parameter Brownian sheet. 

Let us denote by $\mathbb{R}^2_+$ the positive quadrant of the plane and let $z\in \mathbb{R}^2_+$,
%given an $\mathcal{F}_{t,x}$-adapted process $\{\phi(t,x):  t \geq 0, x \in \mathbb{R}\}$. 
We define a first type stochastic integral with respect to the two-parameter Brownian motion in Cairoli \cite{Cairoli72} denoted by:
\begin{align*}
   \int_{R_z} \phi(\zeta)B(d\zeta)
\end{align*}
and a second type \cite{WZ74} stochastic integral denoted by
\begin{align*}
   \int_{R_z} \int_{R_z}\psi(\zeta,\zeta')B(d\zeta)B(d\zeta'),
\end{align*}
where $R_z=[0,t]\times [0,x]$ if $z=(t,x)$.\\

In Wong \& Zakai \cite{WZ}, an It\^{o} formula for stochastic integrals in the plane is given.

\begin{example}{(An optimal harvesting problem)} \\

A classical model for the growth of an ecosystem (e.g. a population or a forest) with value $Y(t)$ at time $t$ in a random environment is the geometric Brownian motion, defined by the It\^{o} stochastic differential equation (SDE)
\begin{equation*}
dY(t)=\kappa Y(t)dt+\gamma Y(t)dB(t), \quad t\geq 0,\quad Y(0) > 0, %\label{(7)}%,
\end{equation*}
where $\kappa$ and $\gamma$ are given constants.
Equivalently, in terms of white noise $\overset{\bullet}{B}$ and Wick product $\diamond$, the equation can be written
\begin{equation*}
\frac{d}{dt}Y(t)=\kappa Y(t)+\gamma Y(t)\diamond\overset{\bullet}{B}%
(t),\quad Y(0)>0,%\label{(8)}%
\end{equation*}
where%
\begin{equation*}
\overset{\bullet}{B}(t)=\frac{d}{dt}B(t)\text{ is (time) white noise.}%\label{(9)}%
\end{equation*}
\end{example}

A natural extension of this model to the case where the noise of the environment depends on both time $t$ and position $x$, is the following SPDE in the value $Y(t,x)$ of the ecosystem at time $t$ and position $x$:
\begin{equation} 
\frac{\partial^{2}}{\partial t\partial x}Y(t,x)=\alpha_{0}(t,x)Y(t,x)+\beta
_{0}(t,x)Y(t,x)\diamond\overset{\bullet}{B}(t,x),\quad Y(0,0) > 0, \label{Y3}%
\end{equation}
where $\alpha_0(t,x)$ and $\beta_0(t,x)$ are given bounded deterministic functions, and
\begin{align} \label{WN}
\overset{\bullet}{B}(t,x)=\frac{\partial^{2}}{\partial t\partial x}B(t,x)\text{ is
time-space white noise.}
\end{align}

Alternatively, if we use the following result:
\small
\begin{equation}
\int_{R_z}\varphi(s,a)B(ds,da)=\int_{R_z}\varphi(s,a
)\diamond\overset{\bullet}{B}(s,a)dsda,\text{ } \quad \forall\varphi, t, x.\label{4a}%
\end{equation}
(see e.g. Holden et al \cite{HOUZ}),
we can rewrite the equation \eqref{Y3} as:
\begin{align}\label{Y4}
Y(t,x)=Y(0,0) + \int_{R_z} \alpha_0(\zeta)Y(\zeta) d\zeta + \int_{R_z} \beta_0(\zeta) Y(\zeta)B(d\zeta).
\end{align}

Assume for simplicity that $\alpha_0$ and $\beta_0$ are constants.\\
If we at $(t,x)$ harvest from $Y(t,x)$ at the rate $u(t,x)$, 
the dynamics (\ref{Y3}) becomes%
\begin{equation*}
\frac{\partial^{2}}{\partial t\partial x}Y_{u}(t,x)=\alpha_{0}Y_{u}%
(t,x)-u(t,x)+\beta_{0}Y_{u}(t,x)\diamond\overset{\bullet}{B}(t,x),%\label{(11)}%
\end{equation*}
or, in integral form,%
\begin{align*}
Y_{u}(t,x)  & = Y(0,0) +\int_0^t\int_0^{x}\{\alpha_{0}Y_u(s,a)-u(s,a)\}dsda \\ %\label{(12)}
&+\int_0^t \int_0^{x}\beta_{0}Y_{u}(s,a)B(ds,da).\nonumber
\end{align*}

For given utility functions $U_1, U_2$ and given constants $T>0, X>0$ 
%such that $T>t, X>x$, 
define the combined utility of the harvesting and the terminal population by
\begin{equation*}
J(u)=E\left[\int_{0}^{T}\int_{0}^{X}U_{1}(u(s,a))dsda+
U_{2}(Y_{u}(T,X))\right].%\label{(13)}%
\end{equation*}
We want to find the harvesting strategy $u^{\ast}(s,x)$ which maximizes the
utility of the harvest, i.e.
\[
J(u^{\ast})=\sup_{u\in\mathcal{A}}J(u).
\]
Here the set of admissible controls is denoted by $\mathcal{A}.$
%\end{example}

\begin{example}{(A linear-quadratic (LQ) problem)}\\
Consider the following linear-quadratic (LQ) control problem for time-space random fields:\\
Suppose the state $Y(t,x)$ is given by%
\begin{align*} %\label{LQ1}
Y(t,x)=Y(0,0)+\int_{0}^t\int_{0}^{x}u(s,a)ds da+\beta B(t,x),\quad t\geq 0, x \in \mathbb{R}.
\end{align*}
We want to drive the state $Y$ to 0 at time-space $(T,X)$ with minimal use of energy. Hence
we put%
\begin{equation*}
J(u)=-\tfrac{1}{2} E\Big[\int_{0}^{T}\int_{0}^{X}u^{2}(s,a)dsda+\theta Y^{2}(T ,X)\Big],%\label{LQ2}%
\end{equation*}
where $\theta >0$ is a given constant.\\
The problem is to find $u^{\ast}\in \mathcal{A}$ such that%
\begin{align}
J(u^{\ast})
=\sup_{u\in\mathcal{A}}J(u).
\end{align}
\end{example}

\begin{example}{(A machine learning problem)}\\
Consider the following hyperbolic SPDE:
\small{
\begin{equation}
Y(t,x)=y-\int_{0}^{t}\int_{0}^{x}u(s,a)\nabla f(Y(s,a))dsda+\sigma
B(t,x),y\in \mathbb{R}^{d},t,x\geq 0\text{,}  \label{HS}
\end{equation}
}
where $B$ is a Brownian sheet in $\mathbb{R}^{d}$, $\sigma \in \mathbb{R}%
^{d\times d}$ and $\nabla f $ is the gradient of a function $f\in C^{1}(\mathbb{R}%
^{d};\mathbb{R})$. Further, $u:\Omega \times \left[ 0,\infty \right)
^{2}\longrightarrow \left[ 0,\infty \right) $ is a stochastic learning rate
in time \emph{and} space, which is assumed to be an adapted random field. By
formally setting $u=\eta \delta _{x}$ in (\ref{HS}) for the Dirac delta
function $\delta _{x}$ in a fixed point $x$ and $\eta \geq 0$ we get
(as a special case of (\ref{HS}) the SDE%
\begin{equation}
dY_{t}=-\eta \nabla f(Y_{t})dt+\sigma dB_{t},Y_{0}=x\in \mathbb{R}^{d},t\geq
0\text{.}  \label{MLS}
\end{equation}
\end{example}

The latter type of SDE (\ref{MLS}) is used in
machine learning in connection with the stochastic gradient descent method
(SGD) to minimize or maximize the objective or loss function $f$. Since the 
dynamics (\ref{HS}) is more general than that of (\ref{MLS}), one may
replace for the sake of a deeper understanding of the classical\ SGD
approach (at the possible expense of numerical tractability) the equation (%
\ref{MLS}) by (\ref{HS}) and study an optimal control problem with respect
to the (time-space) stochastic learning rate $u$. 

In order to illustrate this application
in a simplified framework, we approximate $\nabla f$ \ (for smooth $f$) by
Taylor`s expansion in the case of $d=1$ and consider in this paper the
controlled process
\small{
\begin{align*}
Y_{u}(t,x)&=Y(0,0)-\int_{0}^{t}\int_{0}^{x}u(s,a)Y_{u}(s,a)dsda\\
&+\beta
_{0}B(t,x),Y(0,0),\beta _{0}\in \mathbb{R},t,x\geq 0
\end{align*}
}
with respect to (in this context) natural performance functional%
\begin{equation*}
J(u)=-E\left[ \int_{0}^{T}\int_{0}^{X}u^{2}(s,a)dsda+\theta Y^{2}(T,X)\right]
\end{equation*}
for $\theta >0$.
%\end{example}

\section{A general formulation}

The examples mentioned above are special cases of the following general optimal stochastic control problem:\vskip 0.2cm
We study optimal control of solutions $Y(t,x),t\geq
0,x\in\mathbb{R}$ of SPDEs of the form
\small{
\begin{equation}
Y_{u}(t,x)=Y(t_{0},x_{0})+\int_{R(t,x)}\alpha_{u}(Y_{u}(s,a))dsda
+\int_{R(t,x)}\beta_{u}(Y_{u}(s,a))B(ds,da),\label{(1)}%
\end{equation}
}
where%
\[
R(t,x)=R^{(t_0,x_0)}(t,x)=[t_{0},t] \times[x_{0},x],t\geq t_{0},x\geq x_{0},
\]
and
\[
B(t,x)=\text{Brownian sheet.}
\]

The differential form of (\ref{(1)}) is%
\begin{equation}
\frac{\partial^{2}}{\partial t\partial x}Y_{u}(t,x)=\alpha_{u}(Y_{u}%
(t,x))+\beta_{u}(Y_{u}(t,x))\diamond\overset{\bullet}{B}(t,x).\label{(2)}%
\end{equation}
The identity of (\ref{(1)}) and (\ref{(2)}) comes from the fact that%
\begin{equation}
\int_{R(t,x)}\varphi(s,a)B(ds,da)=\int_{R(t,x)}\varphi(s,a
)\diamond\overset{\bullet}{B}(s,a)dsda,\text{ } \quad \forall\varphi, t, x.\label{(4)}%
\end{equation}
See e.g. Holden et al \cite{HOUZ}.\\

\begin{remark}
Hyperbolic SPDEs of the type \eqref{Y3} have been studied
over the years by several authors. See e.g. Cairoli \cite{Cairoli72} and Yeh 
\cite{Yeh81}, who established strong existence and pathwise uniqueness of
solutions $Y$ of the equation
\small
\begin{align}
Y(t,x)&=y_0+\int_{R(t,x)}b(s,a,Y(s, a))dsda\nonumber\\
&+\int_{R(t,x)} \sigma
(s, a,Y(s,a))B(ds,da), y_0\in \mathbb{R}%
^{d},  \label{HSPDE}
\end{align}%
when $b$ and $\sigma $ are Lipschitz continuous vector fields of linear
growth. 
\end{remark}

Further, smoothness of solutions to (\ref{HSPDE}) in the sense of
Malliavin differentiability for sufficiently regular  $b$ and $\sigma $ was
analysed in Nualart \& Sanz \cite{NualartSanz}. See also Bogso et al \cite{BDMPP}, where the authors
construct Malliavin differentiable unique solutions to (\ref{HSPDE}),
when the drift vector field $b$ is merely bounded and measurable and $\sigma 
$ is given by the unit matrix. As for other works in this direction (in the
case of both weak and strong solutions), we also refer to Yeh \cite{Yeh87} and \cite%
{Yeh85}.\\
\vskip 0.2cm
In addition we refer to the book of Nualart \cite
{Nualart}.\\
 Finally, we want to point out
the interesting link between hyperbolic SPDEs and non-linear (random) wave
equations, when $d=1$ and $b$, $\sigma :\mathbb{R}\mapsto \mathbb{R}$:

By applying the orthogonal transformation $u=x+t, v=x-t$ to the SPDE \eqref{(2)} we see that the corresponding differential
version of equation (\ref{HSPDE}) can be transformed into the following 
%\\ Using formal $\frac{\pi }{2}$ rotation, the corresponding differential
%version of equation (\ref{HSPDE}) can be written as the following
non-linear stochastic wave equation:%
\begin{equation*}
\frac{\partial ^{2}}{\partial t^{2}}Y(t,x)-\frac{\partial ^{2}}{\partial
x^{2}}Y(t,x)=\sigma (Y(t,x))\frac{\partial ^{2}}{\partial t\partial x}
\widetilde{B}(t,x)+b(Y(t,x)),
\end{equation*}
where $\widetilde{B}$ is another Brownian sheet, obtained by applying the inverse orthogonal transformation to $B$.
See e.g. Walsh \cite{Walsh} 
%and Farr\'{e} and Nualart \cite{FarreNualart} 
for
further details.
%\end{remark}

In the sequel we assume that the performance functional $J_{u}(t_{0},x_{0})$ has the form
\begin{equation}
\small{J_{u}(t_{0},x_{0})=E^{(t_{0},x_{0})}\left[\int_{t_0}^T\int_{x_0}^X f_{u}(Y_{u}(s,a
))dsda+g(Y_{u}(T,X))\right],\label{(5)} }
\end{equation}
and consider the following problem:

\begin{problem}\label{Problem}
Let $\mathcal{A}$ be a given family of admissible controls. Find an optimal control $u^{\ast}%
\in\mathcal{A}$ and the value function $\Phi$ such that%
\begin{equation}
\Phi(t_0,x_0)= J_{u^{\ast}}(t_{0},x_{0})=\sup_{u\in\mathcal{A}}J_{u}(t_{0},x_{0}).\label{(6)}%
\end{equation}
The control $u^{\ast}$ is called an optimal control and the function $\Phi$ is referred to as the value function of this problem.
\end{problem}

We handle those problems by using a maximum principle approach. Therefore we need to study the adjoint equation given by a BSPDE in the plane. The existence and the uniqueness of solutions to such a BSPDE were proven in Zaidi \& Nualart \cite{ZN} for a particular Lipschitz constant.\\

Here is an outline of the rest of the paper:
\begin{itemize}
    \item 
In the next section we check properties of the solution of SPDE.
\item In Section 3 we introduce some background about stochastic calculus of time-space white noise. 
\item
In Section 4 we prove two types of maximum principles,  a sufficient (verification theorem) and a necessary maximum principle.
\item
Finally, in Section 5 we apply the results of Section 4 to study the problems mentioned above.
\end{itemize}

%%%%%%%%%%%%%%%%%%%%%%%%%%%%
%%%%%%%%%%%%%%%%%%%%%%%%%%%%%%%%
\section{A  discussion of the solution of the SPDE }
%\frametitle{A  discussion of the solution of the SPDE }
 It is natural to ask what are the properties the solution $Y(t,x)$ of the equation
 \begin{align}\label{Y4a}
Y(t,x)=Y(0,0) + \int_{R_z} \alpha_0(\zeta)Y(\zeta) d\zeta + \int_{R_z} \beta_0(\zeta) Y(\zeta)B(d\zeta).
\end{align}
For example, is $Y(t,x)$ always positive?

 To study this question, let us suppose that the coefficients $\alpha _{0},\beta _{0}$ are bounded
deterministic functions. Then it follows from \eqref{Y3} that%
\begin{equation*}
\widetilde{Y(t,x)}(z)=Y(0,0)+\int_{R_{t,x}}\widetilde{K(s,a)}(z)\widetilde{%
Y(s,v)}(z)dsdv\text{,}
\end{equation*}%
where $\widetilde{Y(t,x)}(z)$ denotes the Hermite transform of $Y(t,x)$ for $%
z\in \left( \mathbb{C}^{\mathbb{N}}\right) _{c}$ (set of all finite
sequences in $\mathbb{C}^{\mathbb{N}}$) and where $R_{(t,x)}=\left[ 0,t\right]
\times \left[ 0,x\right] $. Further, 
\begin{equation*}
K(s,a)=\alpha _{0}(s,a)+\beta _{0}(s,a)\overset{\bullet }{B}(s,a),
\end{equation*}%
where%
\begin{equation}
\overset{\bullet }{B}(s,a)=\frac{\partial ^{2}}{\partial s\partial a}B(s,a)%
\text{ is }\text{white noise.}  \label{(3)}
\end{equation}%
See \cite{HOUZ} for the properties of the Hermite transform.

Then, using Picard iteration, we find with $y_{0}=Y(0,0)>0$, the
semi-explicit solution
\small
\label{Y5}
\begin{align}
&\widetilde{Y(t,x)}(z) \nonumber\\
&=y_{0}\sum_{n=0}^{\infty
}\int_{R_{(t,x)}}\int_{R_{(s_{1},a_{1})}}...\int_{R_{(s_{n-1},a_{n-1})}}\prod%
_{j=1}^{n}\widetilde{K(s_{j},a_{j})}(z)ds_{n}da_{n}...ds_{1}da_{1},  
\end{align}
for $z\in \left( \mathbb{C}^{\mathbb{N}}\right) _{c}$, where $s_{1}>...>s_{n}$, $a_{1}>...>a_{n}$. \\
\vskip 0.3cm

By Theorem 2.11.4 in Holden et al \cite{HOUZ} $Y(t,x)$ is positive if and only if for all $m$ the function 
\begin{equation}
g(y):=\widetilde{Y(t,x)}(iy)e^{-\tfrac{1}{2}|y|^{2}};\quad i=\sqrt{-1}%
,y=(y_{1},y_{2},...,y_{m})\in \mathbb{R}^{m}  \label{PositiveD}
\end{equation}%
is positive definite. 

In this context, note that 
\begin{align}
\widetilde{K(s,a)}(iy)& =\alpha _{0}(s,a)+\beta _{0}(s,a)\widetilde{\overset{%
\bullet }{B}(s,a)}(iy)  \notag \\
& =\alpha _{0}(s,a)+\beta _{0}(s,a)\sum_{k=1}^{m}\mu _{k}(s,a)iy_{k}.
\end{align}%
Here $\mu _{k}(s,a)$ is the $k$'th element in an orthonormal basis of $L^{2}(%
\mathbb{R}^{2})$ consisting of tensor products of Hermite functions. See
Holden et al \cite{HOUZ}, Section 2.2.1 and 2.2.2.

Combining the above we get%
\begin{eqnarray}
&&\widetilde{Y(t,x)}(iy)  \notag \\
&=&y_{0}\sum_{n=0}^{\infty
}\int_{R(t,x)}\int_{R(s_{1},a_{1})}...\int_{R(s_{n-1},a_{n-1})}  \notag \\
&&\prod_{j=1}^{n}\left( \alpha _{0}(s_{j},a_{j})+i\beta
_{0}(s_{j},a_{j})\sum_{k=1}^{m}\mu _{k}(s_{j},a_{j})y_{k}\right)
ds_{n}da_{n}...ds_{1}da_{1}\text{.}
\end{eqnarray}

Therefore the positivity question is equivalent to the following: \newline

\emph{\ Is for all $m=1,2.,...$ the function $g:\mathbb{R}^{m}\mapsto 
\mathbb{C}$ given by}%
\begin{eqnarray*}
&&g(y):=y_{0}\sum_{n=0}^{\infty
}\int_{R(t,x)}\int_{R(s_{1},a_{1})}...\int_{R(s_{n-1},a_{n-1})} \\
&&\prod_{j=1}^{n}\left( \alpha _{0}(s_{j},a_{j})+i\beta
_{0}(s_{j},a_{j})\sum_{k=1}^{m}\mu _{k}(s_{j},a_{j})y_{k}\right)
ds_{n}da_{n}...ds_{1}da_{1}e^{-\tfrac{1}{2}|y|^{2}}
\end{eqnarray*}%
\emph{{\small \ }positive definite}?\emph{\small \ }

It turns out that the latter is not true, in general. In what follows, we
want to give an explanation for this in the case of $\alpha _{0}=0$ and $%
\beta _{0}$ given by%
\begin{equation*}
\beta _{0}(s,a)=\beta _{1}(s)\beta _{2}(a)\text{,}
\end{equation*}%
where $\beta _{1}$ and $\beta _{2}$ are bounded measurable functions. Assume
that $m=1$. We also note that we can write%
\begin{equation*}
\mu _{1}(s,a)=\xi _{1}(s)\xi _{2}(a)
\end{equation*}%
with elements $\xi _{1}$ and $\xi _{2}$ of an orthonormal basis of {\small $%
L^{2}(\mathbb{R})$.} 

In this case, we obtain the representation%
\small
\begin{eqnarray}
&&\widetilde{Y(t,x)}(iy_{1})  \notag \\
&=&y_{0}\sum_{n=0}^{\infty
}\int_{R(t,x)}\int_{R(s_{1},a_{1})}...\int_{R(s_{n-1},a_{n-1})}\prod%
_{j=1}^{n}\left( i\beta _{0}(s_{j},a_{j})\mu
_{1}(s_{j},a_{j})y_{1}\right) ds_{n}da_{n}...ds_{1}da_{1}  \notag \\
&=&y_{0}\sum_{n=0}^{\infty }\frac{1}{(n!)^{2}}\left( i\eta (t,x)y_{1}\right)
^{n}\text{,}  \label{Y2}
\end{eqnarray}%
{\small \ }where{\small \ }$\eta (t,x):=\int_{R(t,x)}\beta _{0}(s,a)\mu
_{1}(s,a)dsda$. Let us now have a look at the expression%
\begin{equation*}
i^{k}\left( \eta (t,x)y_{1}\right) ^{k}\exp (-\frac{1}{2}\left\vert
y_{1}\right\vert ^{2}).
\end{equation*}%

Then for a standard normally distributed random variables $Z$, the inverse
Fourier transform of the latter is given by%
\begin{eqnarray*}
&&i^{k}\frac{1}{(\sqrt{2\pi })^{1/2}}\int_{\mathbb{R}}\left( \eta
(t,x)y_{1}\right) ^{k}\exp (-\frac{1}{2}\left\vert y_{1}\right\vert
^{2})\exp (iu_{1}y_{1})dy \\
&=&i^{k}E\left[ \left( \eta (t,x)Z\right) ^{k}\exp (iu_{1}Z)\right]  \\
&=&\left( \eta (t,x)\right) ^{k}\frac{\partial ^{k}}{\partial u_{1}^{k}}%
\varphi _{Z}(u_{1})=\left( \eta (t,x)\right) ^{k}\frac{\partial ^{k}}{%
\partial u_{1}^{k}}(\exp (-\frac{1}{2}u_{1}^{2})) \\
&=&\left( \eta (t,x)\right) ^{k}(-1)^{k}h_{k}(u_{1})\exp (-\frac{1}{2}%
u_{1}^{2})\text{,}
\end{eqnarray*}%
where $\varphi _{Z}$ denotes the characteristic function of{\small \ }$Z$%
{\small \ }and where{\small \ }$h_{k}$ is the $k$'th Hermite polynomial.%
{\small \ \newline
}So, using dominated convergence, the inverse Fourier transform of the
function $g$ (for $m=1$) is%
\begin{equation}
b(u_{1})\exp (-\frac{1}{2}u_{1}^{2})\text{,}  \label{bexp}
\end{equation}%

where%
\begin{equation*}
b(u_{1}):=y_{0}\sum_{n=0}^{\infty }\frac{1}{(n!)^{2}}\left( \eta
(t,x)\right) ^{n}(-1)^{n}h_{n}(u_{1})\text{.}
\end{equation*}%
Let us now show that the function $b$ cannot be non-negative. For this
purpose let us recall the following properties of Hermite polynomials:%
\begin{equation*}
h_{n}(x+y)=\sum_{k=0}^{n}\frac{n!}{(n-k)!k!}x^{n-k}h_{k}(y)
\end{equation*}%
and%
\begin{equation*}
h_{k}\left( \int_{\mathbb{R}^{2}}\phi (s,a)B(ds,da)\right) =\left( \int_{%
\mathbb{R}^{2}}\phi (s,a)B(ds,da)\right) ^{\diamond k}\text{,}
\end{equation*}%
where $\diamond $ denotes the Wick product and where $\phi \in L^{2}(\mathbb{%
R}^{2})$ with $\left\Vert \phi \right\Vert _{L^{2}(\mathbb{R}^{2})}=1$ (see 
\cite{HOUZ}).

So we obtain that%
\small
\begin{align*}
&h_{n}\left(c+\int_{\mathbb{R}^{2}}\phi (s,a)B(ds,da)\right) \\
&=\sum_{k=0}^{n}\frac{n!}{%
(n-k)!k!}c^{n-k}h_{k}\left(\int_{\mathbb{R}^{2}}\phi (s,a)B(ds,da)\right) \\
&=\left( c+\int_{\mathbb{R}^{2}}\phi (s,a)B(ds,da)\right) ^{\diamond n}
\end{align*}%
for all $c\in \mathbb{R}$. Assume that the above function $b$ is
non-negative. Then the random variable%
\begin{align*}
&b\left(c+\int_{\mathbb{R}^{2}}\phi (s,a)B(ds,da)\right) \\
&=y_{0}\sum_{n=0}^{\infty }%
\frac{(-1)^{n}}{(n!)^{2}}\left( \eta (t,x)\right) ^{n}h_{n}\left(c+\int_{\mathbb{%
R}^{2}}\phi (s,a)B(ds,da)\right) \\
&=y_{0}\sum_{n=0}^{\infty }\frac{1}{(n!)^{2}}\left( -\eta (t,x)(c+\int_{%
\mathbb{R}^{2}}\phi (s,a)B(ds,da))\right) ^{\diamond n}\text{.}
\end{align*}
must be non-negative, too, for all $c\in \mathbb{R}$.

Let $\phi $ be the
first basis element of an orthornormal basis of (tensored) Hermite functions
in $L^{2}(\mathbb{R}^{2})$. Then the Hermite transform of $b(c+\int_{\mathbb{%
R}^{2}}\phi (s,a)B(ds,da))$ evaluated at $z=iy$, $y\in \mathbb{R}$ is given
by%
\begin{equation*}
y_{0}\sum_{n=0}^{\infty }\frac{1}{(n!)^{2}}\left( -\eta (t,x)(c+iy)\right)
^{n}\text{.}
\end{equation*}%
So using once again the positivity criterion in (\ref{PositiveD}), the
expression%
\begin{equation*}
y_{0}\sum_{n=0}^{\infty }\frac{1}{(n!)^{2}}\left( -\eta (t,x)(c+iy)\right)
^{n}e^{-\tfrac{1}{2}|y|^{2}}
\end{equation*}%
must be \emph{in particular} for $y=0$ non-negative (according to the
definition of a positive-definite function). 

So for%
\begin{equation*}
f_0(r):=\sum_{n=0}^{\infty }\frac{1}{(n!)^{2}}r^{n}
\end{equation*}%
the expression%
\begin{equation*}
y_{0}f_0(-\eta (t,x)c)
\end{equation*}%
must be non-negative for all $c$. On the other hand, it is known that the
function $f_0$, which is related to the Bessel function of order zero, has the
following asymptotic property:%
\begin{equation*}
f_0(y)\sim (\pi \sqrt{\left\vert y\right\vert })^{-1/2}\cos (2\sqrt{\left\vert
y\right\vert }-\frac{\pi }{4})\text{ for }y\rightarrow -\infty \text{.}
\end{equation*}%

See e.g. Nualart \cite{Nualart87}.\ Hence, if $\eta (t,x)\neq 0$, the above
expression can become negative, which leads to a contradiction.\ So we
conclude that the function $b$ and therefore also the function in (\ref{bexp}%
) can become negative. So it follows from Bochner's theorem that the
function $g$ is not positive definite.

We get the following conclusion:

\begin{theorem}
The solution $Y(t,x)$ of the equation \eqref{Y5}  is not always positive
\end{theorem}

\begin{remark}
Using a different method, we mention that the latter result was obtained in
Nualart \cite{Nualart87} in the special case of $\alpha _{0}=0$ and $\beta
_{0}$ a constant. There Nualart proves that $Y(t,x)<0$ uniformly in $(t,x)
$ on an open set, with positive probability.
\end{remark}

\section{Background}
%\frametitle{Background}

\subsection{The It\^{o} formula}
To study such optimal control problems, we will use a version of the It\^{o}
formula for such systems. First we introduce some notation
 from Wang \& Zakai \cite{WZ}.
\begin{itemize}
\item
We put $ \zeta=(\zeta_1,\zeta_2)=(s,a)\in \mathbb{R} \times \mathbb{R}$ and 
$d\zeta=d\zeta_{1}d\zeta_{2}=ds da$,
\item
$B(t,x)$ is a Brownian sheet; $t\geq0,x\in\mathbb{R}$,
\item
$z=(z_{1},z_{2})=(t,x),R_{z}=[0,z_{1}]\times\lbrack0,z_{2}]$,
\item
$\int_{R_{z}}\varphi(\zeta)B(d\zeta)$ 
\text{denotes the It\^{o}
integral with respect to} 
$B(\cdot)$ \text{over} $R_{z}$,
\item
$\int_{R_{z}}\psi(\zeta)d\zeta$ is 2-dimensional Lebesgue integral of $\psi$, 

\item
if $a=(a_{1},a_{2}),b=(b_{1},b_{2})$, then
$a\vee b=(\max(a_{1},b_{1}),\max(a_{2},b_{2})).
$
\item
$I((a_1,a_2)\bar{\wedge} (b_1,b_2))=1 \newline
\text{  if } a_1 \leq b_1\quad  \& \quad a_2 \geq b_2,\quad 0 $ \text {otherwise}.
\end{itemize}

\begin{theorem}
[It\^{o} formula, Wang \& Zakai \cite{WZ}]

Suppose%
\begin{equation}
Y(z)=Y_{0}+\int_{R_{z}}\alpha(\zeta)d\zeta+\int_{R_{z}}\beta(\zeta
)B(d\zeta) + \iint\limits_{R_{z}\times R_{z}} \psi (\zeta,\zeta^{\prime})B(d\zeta)B(d\zeta^{\prime}).\label{(16)}%
\end{equation}

\end{theorem}

Then, if $f:\mathbb{R\rightarrow R}$ is smooth, we have

\small
\begin{align*}\label{(17)}
&f(Y(z))  =f(Y_{0})+\int_{R_{z}}f^{\prime}(Y(\zeta))[\alpha(\zeta
)d\zeta+\beta(\zeta)B(d\zeta)]\\
&+\tfrac{1}{2}\int_{R_{z}}f^{\prime\prime}(Y(\zeta))\beta^{2}(\zeta
)d\zeta\nonumber\\
& +\iint\limits_{R_{z}\times R_{z}} \Big\{ f^{\prime\prime}(Y(\zeta\vee\zeta^{\prime
}))u_k \tilde{u}_l+f^{\prime} (Y(\zeta\vee \zeta^{\prime})) \psi(\zeta,\zeta^{\prime}) \Big\}B(d\zeta)B(d\zeta^{\prime}) \nonumber\\
&+\iint\limits_{R_{z}\times R_{z}}\Big\{f^{\prime\prime}(Y(\zeta\vee\zeta
^{\prime}))\Big( u_k \alpha(\zeta) + \psi(\zeta,\zeta^{\prime})  \tilde{u}_l \Big)\nonumber\\
& +\tfrac{1}{2}f^{(3)}(Y(\zeta\vee\zeta^{\prime}))u^2_k \tilde{u}_l\Big\}d\zeta B(d\zeta^{\prime})\nonumber\\
&+\iint\limits_{R_{z}\times R_{z}}\Big\{f^{\prime\prime}(Y(\zeta\vee\zeta
^{\prime}))\Big( \tilde{u}_l\alpha(\zeta^{\prime}) + \psi(\zeta,\zeta^{\prime}) u_k  \Big)\nonumber\\
& +\tfrac{1}{2}f^{(3)}(Y(\zeta\vee\zeta^{\prime}))u^2_k \tilde{u}_l\Big\}d\zeta B(d\zeta^{\prime})\nonumber\\
& +\iint\limits_{R_{z}\times R_{z}}I(\zeta \bar{\wedge} \zeta^{\prime}) \Big\{f^{\prime\prime}(Y(\zeta\vee\zeta
^{\prime}))\Big(\alpha(\zeta^{\prime})\alpha(\zeta)+\tfrac{1}{2}\psi^2(\zeta,\zeta^{\prime})\Big) \nonumber\\
&+f^{(3)}(Y(\zeta\vee\zeta^{\prime}))u \tilde{u}\psi(\zeta,\zeta^{\prime})+\tfrac{1}{2}f^{(3)}(Y(\zeta\vee\zeta^{\prime}))\left[\alpha(\zeta^{\prime}%
)\tilde{u}^2+\alpha(\zeta)u^2 \right]\nonumber\\
& + \tfrac{1}{4}f^{(4)}(Y(\zeta\vee\zeta^{\prime}))u^2 \tilde{u}^2\Big\}d\zeta d\zeta^{\prime},\nonumber
\end{align*}
%\end{theorem}

where

$$u=\beta(\zeta^{\prime})+\int_{R_{z}}I(\zeta \bar{\wedge} \zeta^{\prime})\psi(\zeta,\zeta^{\prime})B(d\zeta),$$

$$\tilde{u}=\beta(\zeta)+\int_{R_{z}}I(\zeta \bar{\wedge} \zeta^{\prime})\psi(\zeta,\zeta^{\prime})B(d\zeta^{\prime}).$$
%\end{theorem}

It is proved in \cite{WZ} that the double $B(d\zeta)B(d\zeta')$-integrals, and the mixed $d\zeta B(d\zeta')$ and $B(d\zeta)d\zeta'$-integrals are all weak martingales and hence have expectation $0$. Therefore, by the It\^{o} formula above we get the following:
\begin{theorem}
(Dynkin formula)
\small
\begin{align*}
E[f(Y(z))]  & =f(Y_{0})+E\Big[\int_{R_{z}}\Big\{\alpha(\zeta)f^{\prime}(Y(\zeta
))+\frac{1}{2}\beta^{2}(\zeta)f^{^{\prime\prime}}(Y(\zeta))\Big\}d\zeta\\
&+\iint\limits_{R_{z}\times R_{z}}I(\zeta \bar{\wedge} \zeta^{\prime}) \Big\{f^{\prime\prime}(Y(\zeta\vee\zeta
^{\prime}))\Big(\alpha(\zeta^{\prime})\alpha(\zeta)+\tfrac{1}{2}\psi^2(\zeta,\zeta^{\prime})\Big) \\
&+f^{(3)}(Y(\zeta\vee\zeta^{\prime}))u_k \tilde{u}_l\psi(\zeta,\zeta^{\prime})\\
&+\tfrac{1}{2}f^{(3)}(Y(\zeta\vee\zeta^{\prime}))\left[\alpha(\zeta^{\prime}%
)\tilde{u}^2_l+\alpha(\zeta)u^2_k \right]\\
& + \tfrac{1}{4}f^{(4)}(Y(\zeta\vee\zeta^{\prime}))u^2_k \tilde{u}^2_l\Big\}d\zeta d\zeta^{\prime}\Big].
\end{align*}
\end{theorem}

\begin{lemma}[Integration by parts]\label{partsa}
Suppose that 
\small
\begin{align*}
    Y_k(z)=Y_k(0)+\int_{R_{z}}\alpha_k(\zeta)d\zeta
    +\int_{R_{z}}\beta_k(\zeta)B(d\zeta)+\iint\limits_{R_{z}\times R_{z}} \psi_k (\zeta,\zeta^{\prime})B(d\zeta)B(d\zeta^{\prime}).\quad k=1,2.
\end{align*}
Then 
\small
\begin{align*}
  &  E[Y_1(z)Y_2(z)]=Y_1(0)Y_2(0)\\
  &+E\Big[\int_{R_{z}}\Big\{ Y_1(\zeta) \alpha_2(\zeta)+Y_2(\zeta) \alpha_1(\zeta)+ \beta_1(\zeta)\beta_2(\zeta)\Big\}d\zeta\nonumber\\
& +\iint\limits_{R_{z}\times R_{z}}I(\zeta \bar{\wedge} \zeta^{\prime}) \Big\{\alpha_1 (\zeta^{\prime}) \alpha_2(\zeta) + \alpha_{1}(\zeta)\alpha_2(\zeta^{\prime}) +\psi_1(\zeta,\zeta^{\prime}) \psi_2(\zeta,\zeta^{\prime}) \Big\} d\zeta d\zeta^{\prime}\Big].\\
\end{align*}
\end{lemma}
\dproof
The proof follows from the It\^o formula, Proposition 5.1 in Wang \& Zakai \cite{WZ}.
\fproof\\

%%%%%%%%%%%%%%%%%%%%%%%%%%%%%%%%%%%%%%%%%%%%%%5%%%%%%%%%%%%%%%%%%%
\section{BSPDEs in the plane}
%\frametitle{BSPDEs in the plane}
To simplify the notation we sometimes put $z=(t,x), \zeta=(s,a)$ in the following:\\
Let $\mathcal{P}$ be the predictable $\sigma$-algebra of subsets of $\Omega \times R_{z_0}$ generated by the sets $(z,z']\times A$ where $A\in \mathcal{F}_z$ and we denote by $\mathcal{D}$ the $\sigma$-algebra of $\Omega \times R_{z_0} \times R_{z_0}$ generated by the sets $(z_1,z_1']\times(z_2,z_2'] \times A$ where $(z_1,z_1']\bar{\wedge}
 (z_2,z_2']$ and $A\in \mathcal{F}_{z_1\vee z_2}$.\\
The solution of the BSDE will live in the following spaces:
\begin{itemize}
    \item $L^2_{a,1}$ is the space of predictable processes \{$\phi(z),z\in R_{z_0}$\}, such that \\$E[\int_{R_z}\phi(z)^2dz<\infty]$,
    \item  $L^2_{a,2}$ is the space of processes \{$\psi(z,z'),(z,z')\in R_{z_0}\times R_{z_0}$\}, such that
    \begin{itemize}
        \item [(a)] $\psi(z,z')=0$ unless $z\bar{\wedge}z'$,
        \item [(b)] $\psi$ is $\mathcal{D}$-measurable,
        \item [(c)]  $E[\int_{R_z} \int_{R_z}\psi(z,z')^2dzdz'<\infty]$.
    \end{itemize}
   
\end{itemize}
Let us recall now the representation of square integrable martingales. 
\begin{theorem}[Wong \& Zakai \cite{WZ}]
  If $M=\{M(z), \mathcal{F}_z, z\in \mathbb{R}_+ ^2\}$ is square integrable martingale, then for each $z\in \mathbb{R}_+ ^2$
  \begin{align}
  M(z)=M(0)+\int_{R_{z}}\phi(\zeta) B(d\zeta)+\iint\limits_{R_{z}\times R_{z}}\psi(\zeta,\zeta') B(d\zeta) B(d\zeta'),
  \end{align}
where $\phi, \psi$ are adapted processes. 
\end{theorem}

Let $Z=(T,X)$ and if we fix a rectangle $R_{Z}=[0,T]\times[0,X]$, and let $\xi$ be an $\mathcal{F}_{Z}$-measurable random variable and $h(\omega,\zeta,p,q)$ is a $\mathcal{P}\times\mathcal{B}_\mathbb{R}\times\mathcal{B}_\mathbb{R}$-measurable function such that $\int_{R_{z}}|h(\zeta,p(\zeta), q(\zeta))|d\zeta<\infty$. Then we can define a triple of processes $(p,q,r)\in L^2_{a,1}\times L^2_{a,1}\times L^2_{a,2}$ solution of the BSPDE in the plane
\begin{align}\label{BSDE1}
p(z)&=\xi-\int_{R_{Z}\setminus R_{z}} h(\zeta,p(\zeta), q(\zeta))d\zeta-\int_{R_{Z}\setminus R_{z}} q(\zeta) B(d\zeta)\nonumber\\
&-\int_{R_{Z}\setminus R_{z}}\int_{R_{Z}\setminus R_{z}} r(\zeta,\zeta') B(d\zeta)B(d\zeta').
\end{align}

Alternatively, let us introduce the notation
\begin{align}
    M_{r}(z)=\iint\limits_{R_{z}\times R_{z}} r(\zeta,\zeta') B(d\zeta)B(d\zeta'),\quad r \in L_{a,2}.
\end{align}
Then $M_{r}(z)$ is a martingale, and we can write the equation for $(p,q,r)$ above in differential form as follows
\begin{align}
        p(dz)&=h(z,p(z), q(z))dz+q(z) B(dz)+M_{r}(dz),\quad  z \leq Z,\label{BSDE2}\\
    p(Z)&=\xi.\nonumber
\end{align}

\textbf{Assumptions} We impose the following set of assumptions:
 \begin{itemize}
        \item [(i)] $\xi \in L^2(\Omega,\mathcal{F}_{z_0},P)$,
       \item [(ii)]$h(\cdot,p,q)\in L^2_{a,1}$ for all $p,q\in \mathbb{R}$, 
        \item [(ii)] $|h(\zeta,p,q)-h(\zeta',p',q')|^2\leq K_1|p-p'|^2+K_2 |q-q'|^2$, for all $p,q,p',q'\in \mathbb{R}$ and  $\zeta \in R_{Z}.$
    \end{itemize}
Let $f_0$ be the Bessel function of order zero and $r_0\approx 1.4458$ be the first nonnegative zero of $J_0$:
$$
r_0=inf \left\{t>0:f_0(2\sqrt{t})=\sum_{j=0}^\infty \frac{(-1)^j}{{j!}^2} t^j=0 \right\}.
$$
\begin{theorem}[Existence and Uniqueness (Zaidi \& Nualart) \cite{ZN}]
Under the above assumptions (i)-(iii) and if the Lipschitz constant satisfies $K_1|z_0|<\sqrt{r_0}$ and $K_2|z_0|<\sqrt{r_0}$, 
there exists a unique solution of the BSPDE \eqref{BSDE1}.
   \end{theorem}

\section{Maximum principle approaches}
%\frametitle{Maximum principles}

We  introduce the following notation:
\begin{definition}
For general functions $h,k:[0,T]\times[0,X] \mapsto \mathbb{R}$, we define
\begin{align}
 (h\star k)(t,x)=\int_0^x \int_{-t}^T \{ h(t,x)k(s,a) + h(s,a) k(t,x)\} ds da  
\end{align}
\end{definition}
Note that with this notation we have
\small
\begin{align}
 \iint\limits_{R_{z}\times R_{z}}I(\zeta \bar{\wedge} \zeta^{\prime}) \Big\{\alpha_1 (\zeta^{\prime}) \alpha_2(\zeta) + \alpha_{1}(\zeta)\alpha_2(\zeta^{\prime})  \Big\} d\zeta d\zeta^{\prime}=\int_{R_{z}}( \alpha_1 \star \alpha_2)(\zeta) d\zeta,   
\end{align}
and the integration by parts formula can be written

\begin{lemma}[Integration by parts (2)]\label{parts}
Suppose that 
\begin{align*}
    Y_k(z)=Y_k(0)+\int_{R_{z}}\alpha_k(\zeta)d\zeta
    +\int_{R_{z}}\beta_k(\zeta)B(d\zeta), \quad k=1,2.
\end{align*}
Then 
\small
\begin{align} 
    &E[Y_1(z)Y_2(z)]=Y_1(0)Y_2(0)\nonumber\\
    &+E\Big[\int_{R_{z}}\Big\{ Y_1(\zeta) \alpha_2(\zeta)+Y_2(\zeta) \alpha_1(\zeta)+ \beta_1(\zeta)\beta_2(\zeta)+(\alpha_1 \star \alpha_2)(\zeta) \Big\} d\zeta \Big].\label{ibp2}
\end{align}
\end{lemma}

Given a subset $U$ of $\mathbb{R}$ we denote by $\mathcal{U}$ the set of all $\mathcal{F}_{t,x}$-adapted control processes $u=\{u(t,x), t<T, x<X\}$ valued in $U$. We then define the set of admissible control processes $\mathcal{A}\subset\mathcal{U}$ to be the collection of all $\mathcal{F}_{t,x}$-adapted processes with values in $U$.\\
 
Let $f$ and $g$ be given functions and consider 
the performance functional
\begin{align*}
J(u)=E\Big[\int_{R_{Z}}f(\zeta,Y(\zeta),u(\zeta))d\zeta+g(Y(Z))\Big],
\end{align*}
where
$$R_Z=[0,T] \times [0,X],$$
with $Z=(T,X)$ for some given $T>0,X>0$, and the state $Y$ of the system is described by the equation
\small
\begin{align}
Y(z)=Y(t,x)=Y(0)+\int_{R_{z}}\alpha(\zeta,Y(\zeta),u(\zeta))d\zeta+\int_{R_{z}}\beta(\zeta,Y(\zeta),u(\zeta))B(d\zeta), \quad z\leq Z,
\end{align}
where $R_z=[0,t] \times [0,x]$ when $z=(t,x)$, and $u$ denotes a control process.

\begin{problem}
We want to find $\widehat{u}\in\mathcal{A}$ such that
\begin{equation}\label{J}
J(\widehat{u})=\sup_{u\in\mathcal{A}}J(u).
\end{equation}
\end{problem}
 The maximum principle approach to this problem, adapted to the time-space situation, is to introduce the following associated Hamiltonian 
$$H: R_{Z} \times R_{Z}\times \R \times U\times \R \times \R \times \R $$
given by
\begin{align}
H(z,Y,u,p,q,L)
=f(z,y,u)+p \alpha(z,y,u)+q \beta(z,y,u)+ (L \star \alpha)(z),
\end{align}
where we recall that  $$ 
(L \star \alpha)(z) =\int_0^x\int_t^T \Big\{L(z) \alpha(\zeta',y,u) + L(\zeta')\alpha(z,y,u)\Big\}d\zeta' ;\quad z=(t,x)$$
and the adjoint processes $(p,q,r,L)=(p(t,x),q(t,x),r(t,x,\cdot),L(t,x))$ are given by the following equations:
 \begin{itemize}
    \item
$L:R_Z \mapsto \mathbb{R}$ is defined implicitly by the integral equation
\begin{align}
  L(z)= -\frac{\partial H}{\partial y} (z,Y(z),u(z), p(z),q(z),r(z,\cdot),L(z)); \quad z=(t,x)
\end{align}
\item
and $(p,q,r)$ is the solution of the BSPDE
\begin{align}
\begin{cases}
    p(dz)&=- \frac{\partial H}{\partial y}(z,Y(z),p(z),q(z),r(z,\cdot),L(z))dz\\
    &+q(z) B(dz) 
     +M_{r}(dz),\quad 0 \leq (t,x) \leq(T,X),\\
    p(Z)&=\frac{\partial g}{\partial y}(Y(Z)),
    \end{cases}
    \end{align}
    \end{itemize}
    
   or, in integrated form,
    \begin{align}
    \begin{cases}
   p(z)&=p(0) - \int_{R_{z}} \frac{\partial H}{\partial y} (\zeta,Y(\zeta),p(\zeta),q(\zeta),r(\zeta,\cdot),L(\zeta)) d\zeta \\
   &+\int_{R_{z}} q(\zeta)B(d\zeta)
   +\iint\limits_{R_{z}\times R_{z}}r(\zeta,\zeta')B(d\zeta)B(d\zeta'),\quad z\leq Z.\\
   p(Z)&=\frac{\partial g}{\partial y}(Y(Z))
   \end{cases}
    \end{align}
   There are two versions of the maximum principle for this problem, namely the so-called \emph{sufficient maximum principle} and the \emph{necessary maximum principle}.
We present them both below.

\subsection{The sufficient maximum principle}
%\frametitle{The sufficient maximum principle}
\begin{theorem}[Sufficient maximum principle]
 Suppose $\widehat{u}\in\mathcal{A}$ with corresponding solutions $\widehat{Y},(\widehat{p},\widehat{q},\widehat{r},\widehat{L})$ of the equations above. Moreover, suppose that
 $
 y 	\mapsto g(y)
 $ is concave and $
 Y,u 	\mapsto H(z,Y, u,p,q,r,L)
 $ is concave for all $z,p,q,r,L$ and that 
 \begin{align}  \label{max.cond}
 \small
     \sup_{v\in A}&H(z,\hat{Y}(\cdot),v,\widehat{p}(z),\widehat{q}(z),\widehat{r}(z,\cdot),\widehat{L}(z))\\
     &=H(z,\hat{Y}(\cdot),\widehat{u}(z),\widehat{p}.(z),\widehat{q}(z),\widehat{r}(z,\cdot),\widehat{L}(z)),
 \end{align}
 Then $\widehat{u}$ is an optimal control for problem \eqref{J}.
\end{theorem}
 
\dproof
Suppose $\widehat{u}\in\mathcal{A}$ satisfies \eqref{max.cond} with corresponding $\widehat{Y}$. Choose another $u\in\mathcal{A}$. Then 
\begin{align}
    J(u)-J(\widehat{u})=I_1+I_2,
\end{align}
where 
\small
\begin{align*}
    I_1&=E\Big[\int_{R_{Z}}\Big\{f(\zeta,Y(\zeta),u(\zeta))-f(\zeta,\widehat{Y}(\zeta),\widehat{u}(\zeta))\Big\}d\zeta\Big]
    =E\Big[\int_{R_{Z}}\Tilde{f}(\zeta)d\zeta\Big],\\
    I_2&=E[g(Y(Z))-g(\widehat{Y}(Z))].
\end{align*}
Using the definition of $H$ we can write 
\begin{align} \label{4.7}
    I_1=E\Big[\int_{R_{Z}}\Big\{ H(\zeta)-\widehat{H}(\zeta)-\widehat{p}(\zeta) \tilde{\alpha}(\zeta)-\widehat{q}(\zeta)\tilde{\beta}(\zeta) -(\widehat{L} \star \Tilde{\alpha})(\zeta)\Big\} d\zeta \Big],
\end{align}
where $\Tilde{\alpha}=\alpha-\widehat{\alpha}, \widehat{\alpha}=\alpha(\zeta,\widehat{Y}(\zeta),\widehat{u}(\zeta))$ etc.

Using the concavity of $g$ and Lemma \ref{parts}, and the fact that the $B(dz)$-integrals and the $B(dz)B(dz')$-integrals are orthogonal (see \cite{CW}, Theorem 2.5), we get 
\begin{align}
    I_2 &\leq E\Big[\frac{\partial g}{\partial y} (\widehat{Y}(Z))\tilde{Y}(Z)\Big]=E\Big[\widehat{p}(Z)\Tilde{Y}(Z)\Big]\nonumber\\
    &=E\Big[\int_{R_{Z}}\Big\{\widehat{p}(\zeta)\tilde{\alpha}(\zeta)
    - \frac{\partial \hat{H}}{\partial y}\Tilde{Y} (\zeta)+\widehat{q}(\zeta)\Tilde{\beta}(\zeta) -(\frac{\partial \hat{H}}{\partial y}\star  \tilde{\alpha})(\zeta)\Big\}d\zeta\Big].\label{4.8}
\end{align} 
Adding \eqref{4.7} and \eqref{4.8} we get, using that
$\hat{L}= - \frac{\partial \hat{H}}{\partial y}$ and the concavity of $H(y,u)$,
\begin{align}
    J(u)-J(\widehat{u}) &=I_1 + I_2 \leq E\Big[\int_{R_{Z}}\Big\{ H(\zeta)-\widehat{H}(\zeta) -\frac{\partial \hat{H}}{\partial y}(\zeta)\Tilde{Y}(\zeta)\Big\}d\zeta\Big] \\ 
    &\leq E\Big[\int_{R_{Z}} \frac{\partial \hat{H}}{\partial u} (\zeta)\Tilde{u}(\zeta)\>d\zeta\Big] \leq 0\text{ by condition \eqref{max.cond}}.
\end{align}
This proves that 
\begin{align}
    J(u)-J(\widehat{u})\leq 0 \text{ for all } u\in\mathcal{A},
\end{align}
and therefore $\widehat{u}$ is optimal.
\fproof

\subsection{The necessary maximum principle}
%\frametitle{The necessary maximum principle}
It is a drawback of the sufficient maximum principle that we have to assume that $
 y 	\longmapsto g(y)
 $ and $
 (y,u) 	\longmapsto H(z,y,u,p,q,L)
 $ are concave. The following result does not need concavity, but we have to add conditions of the set $\mathcal{A}$ of admissible controls instead, as follows: 
\begin{itemize}
    \item [(A1)] $\mathcal{A}$ is a convex set
    \item  [(A2)]  For all $z_0=(t_0,x_0)<Z=(T,X)$ 
    and all bounded $\mathcal{F}_{z_{0}}$-measurable random variables $\theta_{z_{0}}$, the control $$ u_{z_{0}}(\zeta)=\theta_{z_{0}} \mathbf{1}_{R_{z_{0}}}(\zeta)$$
    is admissible, where
    \begin{align}
    \mathbf{1}_{R_{z_{0}}}(\zeta) =
        \begin{cases}
         1 \text{ if } \zeta \in  R_{z_{0}}\\
         0 \text{ if } \zeta \notin  R_{z_{0}}
        \end{cases}
    \end{align}
\end{itemize}
   is the indicator function of the rectangle $R_{z_{0}}=[t_0,T] \times [x_0,X]$.
    
\begin{lemma}
    For all  $u,v \in \mathcal{A}$ the derivative process
    $$G(\zeta):=\underset{\epsilon \rightarrow 0}{\lim} \frac{Y^{u+\epsilon v}(\zeta)-Y^u(\zeta)}{\epsilon}$$
satisfies the equation

\begin{align*}
     G(z)&=G(0)+\int_{R_{z_{0}}} \Big\{ \frac{\partial \alpha}{\partial y} (\zeta) G(\zeta)+ \frac{\partial \alpha}{\partial u} (\zeta) v(\zeta)\Big\} d\zeta \\
     &+\int_{R_{z_{0}}} \Big\{ \frac{\partial \beta}{\partial y} (\zeta) G(\zeta)+ \frac{\partial \beta}{\partial u} (\zeta) v(\zeta)\Big\} B(d\zeta),
 \end{align*}
 where $ \frac{\partial \alpha}{\partial y} (\zeta) =\frac{\partial \alpha}{\partial y} (\zeta,Y^u (\zeta), u(\zeta)) $ etc.
 \end{lemma}

 \dproof
 This follows by the chain rule.
 \fproof
\begin{lemma}
    For all $u,v\in \mathcal{A}$, we have 
   $$  \frac{d}{d\epsilon} J(u+\epsilon v)_{\epsilon=0}=E\Big[\int_{R_{Z}}  \frac{\partial H}{\partial u} (\zeta)v(\zeta)  d\zeta\Big].$$
\end{lemma}
 
\dproof
\small
\begin{align*}
    \frac{d}{d\epsilon} &J(u+\epsilon v)_{\epsilon=0}\\
    &=\underset{\epsilon \rightarrow 0}{lim} \frac{1}{\epsilon} E\Big[\int_{R_{Z}} \Big\{f(\zeta,Y^{u+\epsilon v}(\zeta),u+\epsilon v)(\zeta))-f(\zeta,Y^{u}(\zeta),u)(\zeta))\Big\}d\zeta\nonumber\\
    &+g(Y^{u+\epsilon v}(Z))-g(Y^{u}(Z)) \Big]\\
    &=E\Big[\int_{R_{Z}} \Big\{\frac{\partial f}{\partial y} (\zeta,Y^{u}(\zeta),u(\zeta)) G(\zeta)+\frac{\partial f}{\partial u} (\zeta,Y^{u}(\zeta),u(\zeta)) v(\zeta)\Big\}d\zeta \\
    &+\frac{\partial g}{\partial y} (Y^{u}(Z)) G(Z)\Big]=I_1+I_2,
\end{align*}
 
where
\small
\begin{align*}
  I_1&=E\Big[\int_{R_{Z}} \Big( \frac{\partial H}{\partial y}(\zeta) -\frac{\partial \alpha}{\partial y} (\zeta)p(\zeta)-\frac{\partial \beta}{\partial y} (\zeta) q(\zeta)-(L \star  \frac{\partial \alpha}{\partial y})(\zeta)\Big)G(\zeta) d\zeta 
\\
     &+\int_{R_{Z}}  \Big( \frac{\partial H}{\partial u} (\zeta) -\frac{\partial \alpha}{\partial u} (\zeta)p(\zeta)-\frac{\partial \beta}{\partial u} (\zeta)  q(\zeta)-(L \star \frac{\partial \alpha}{\partial u}(\zeta)\Big) v(\zeta)  d\zeta\Big],
\end{align*}
and
\begin{align*}
I_2&=E\Big[\frac{\partial g}{\partial y} (Y^{u}(Z)) G(Z)\Big]=E[p (Z) G(Z)]\\
&=E\Big[\int_{R_{Z}}\Big(p(\zeta)\Big\{ \frac{\partial \alpha}{\partial y} (\zeta)G(\zeta)
+\frac{\partial \alpha}{\partial u}(\zeta)v(\zeta)\Big\}-\frac{\partial H}{\partial y} (\zeta)G(\zeta)\\
&+ q(\zeta) \Big\{ \frac{\partial \beta}{\partial y}(\zeta) G(\zeta)+\frac{\partial \beta}{\partial u}(\zeta)v(\zeta)\Big\}-\big\{( \frac{\partial H}{\partial y} \star (\frac{\partial \alpha}{\partial y} G + \frac{\partial \alpha}{\partial u}  v)  \big\} (\zeta)\Big\}d\zeta\Big].
\end{align*}
 
Adding $I_1$ and $I_2 $, we get, if we choose
$L= - \frac{\partial H}{\partial y}$,
$$
    \frac{d}{d\epsilon} J(u+\epsilon v)_{\epsilon=0}=E\Big[\int_{R_{Z}} \frac{\partial H}{\partial u}(\zeta) v(\zeta)  d\zeta\Big].
$$
\fproof

From Lemma 6.4, we deduce the following:
\begin{theorem} [Necessary maximum principle]
 Suppose $\widehat{u}\in\mathcal{A}$ is optimal for Problem 2.5. Then 
\begin{align*}  
\frac{\partial H}{\partial u} (\zeta,\widehat{Y}(\zeta),\widehat{u}(\zeta),\widehat{p}(\zeta),\widehat{q}(\zeta),\widehat{L}(\zeta))=0 \text{ for a.a. } \zeta.
\end{align*}  
\end{theorem}

\dproof Since $J(\widehat{u}+\epsilon v)_{\epsilon=0} \leq J(\widehat{u})$ for all $\epsilon,v,$ we get by Lemma 6.4 that 
$$
   E\Big[\int_{R_{Z}} \frac{\partial H}{\partial u}(\zeta)v(\zeta) d\zeta\Big]\leq 0, \text{ for all  } v\in \mathcal{A}.
$$
In particular, applying this to 
$$ v(\zeta)=\theta_{z_{0}} \mathbf{1}_{R_{z_{0}}}(\zeta)$$
as in A2, this gives
$$
   E\Big[\int_{R_{z_0}} \frac{\partial H}{\partial u}(\zeta)\theta_{z_{0}} d\zeta\Big]\leq 0.
$$
Since this holds for all $z_0$ we deduce that
$$
\frac{\partial ^2}{\partial t_0 \partial x_0}  \left( E\left[\int_{R_{z_0}} \frac{\partial H}{\partial u}\theta_{z_{0}} d\zeta\right]\right)=\frac{\partial H}{\partial u}(z_{0}) \theta_{z_{0}}
\leq 0.
$$
Since this holds for all $\theta_{z_{0}} \in\mathbb{R}$, we conclude that
$$
\frac{\partial H}{\partial u}(z_{0})=0.
$$
\fproof

%%%%%%%%%%%%%%%%%%%%%%%%%%%%%%%%%%%%%%%%%%%%%%%%%%%%%%%%%%%%%%%%%%%%%%%%%%%%%%%%%%%%%%%%%%%%%%%%%%%

\section{Applications}
%\frametitle{Applications}
\subsection{Return to the optimal harvesting problem in the plane}

Suppose that the growth of a population at time $t$ and position $x$ with density $Y(t,x)$ satisfies
\small
\begin{align*}
Y_{u}(t,x)  & =Y(0,0)+\int_0^t\int_0^{x}\{\alpha_{0}Y_u(s,a)-u(s,a)\}dsda \\
&+\int_0^t \int_0^{x}\beta_{0}Y_{u}(s,a)B(ds,da),\nonumber
\end{align*}
where $\alpha_{0}, \beta_{0}$ are given constants and $Y(0,0)>0.$

For given constants $T>0, X>0$ such that  $T>t, X>x$, define the combined utility of the harvesting and the terminal population by
\begin{equation*}
J(u)=E\left[\int_{0}^{T}\int_{0}^{X}ln(u^2(s,a))dsda+
\theta Y_{u}(T,X)\right],
\end{equation*}
where $\theta$ is a given bounded, $\mathcal{F}_{Z}$-measurable random variable.\\
 
\begin{problem}\label{7.1}
We want to find the harvesting strategy $u^{\ast}(s,x)$ which maximizes the
utility of the harvest, i.e. we want to find $u^{\ast}\in \mathcal{A}$ such that%
\begin{align*}
J(u^{\ast})
=\sup_{u\in\mathcal{A}}J(u).
\end{align*}
\end{problem}

The associated Hamiltonian for this case is
\begin{align*}
    H(t,x,y,u,p,q,L)=ln(u^2) +(\alpha_0y-u)p +\beta_0yq
    +(L \star (\alpha_0 y-u)).
\end{align*}
The Hamiltonian has a maximum at $u$ given by the equation 
$$\frac{\partial H}{\partial u} = \frac{2}{ u(z)}- p-(L \star 1)(z)=0.$$

Hence we must have
$$
u=\frac{2}{p+(L\star 1)},
$$
where 
\begin{align}
 (L\star 1)(z)=x(T-t)L(z) + \int_0^x\int_t^T L(\zeta')d\zeta' .\label{L1a}
 \end{align}

Since
\begin{align}
     \frac{\partial H}{\partial y}=\alpha_0 p +\beta_0 q +(L \star \alpha_0)
 \end{align} 
 the corresponding adjoint equation is
\begin{align}\label{p1a}
\begin{cases}
    p(dz)&=- (\alpha_0 p +\beta_0 q +L \star \alpha_0)(z)dz\\
    &+q(z) B(dz) 
     +M_{r}(dz),\quad 0 \leq (t,x) \leq(T,X), \\
    p(Z)&=\frac{\partial g}{\partial y}(Y(Z)).
    \end{cases}
    \end{align}
  
The adjoint variable $L$ is defined by
\begin{align}\label{L1b}
    L(z)&=-\frac{\partial H}{\partial y}(z) = -[\alpha_0 p(z) +\beta_0 q(z) +(L\star 1)(z)]\nonumber\\
    &= -[\alpha_0 p(z) +\beta_0 q(z) +x(T-t)L(z)+\int_0^x\int_t^TL(\zeta') d\zeta',
\end{align}
which is an integral equation for $L(z)$.
\vskip 0.3cm
We summarise this as follows:
\begin{theorem}
Let $z>0$ and assume that $p(z)+(L \star 1)(z) \neq 0$. Then the optimal harvesting rate $u^{*}$ for Problem \ref{7.1} is given by
    \begin{align*}
u^{*}(t,x)=u^{*}(z)= \frac{2}{p(z)+(L \star 1)(z)},
    \end{align*}
    with $p$ and $L$ given by \eqref{p1a} and \eqref{L1b}, respectively.
\end{theorem}
 
%\begin{remark}
%If $a_{0}=0$, then using the
%martingale property we get that 
%\begin{equation*}
%E\left[ \frac{\Gamma (Z)}{\theta }\left\vert %\mathcal{F}_{z}\right. \right] =%
%\frac{\Gamma (z)}{\theta }\text{.}
%\end{equation*}%
%We know, however, from Remark \ref{5.3} that $\Gamma (z)$ %has a probability density. So in this case, we see that %$u^{*}(z)=2\theta $ a.e.  
%\end{remark}

\subsection{Return to the linear-quadratic (LQ) problem in the plane}

 We now return to the linear-quadratic (LQ) control problem for time-space random fields discussed in the introduction:
 \vskip 0.2cm
Suppose the state $Y(t,x)$ is given by%
\begin{align} \label{LQ1}
Y(t,x)=Y(0,0)+\int_{0}^t\int_{0}^{x}u(s,a)ds da+\beta B(t,x);\quad t\geq 0, x \in \mathbb{R}.
\end{align}
We want to drive the state $Y(t,x)$ to 0 at time $T$ and point $X$ with minimal use of energy. Hence
we put%
\begin{equation}
J(u)=-\tfrac{1}{2} E\left[\int_{0}^{T}\int_{0}^{X}u^{2}(s,a)dsda+\theta Y^{2}(T ,X)\right],\label{LQ2}%
\end{equation}
where $\theta >0$ is a given constant.
\begin{problem} \label{4.3}
We want to find $u^{\ast}\in \mathcal{A}$ such that%
\begin{align}
J(u^{\ast})
=\sup_{u\in\mathcal{A}}J(u).
\end{align}
\end{problem}
 
In this case we have $\nabla_y  H =0$, so we can simplify the Hamiltonian to
\begin{align}
    H_0(t,x,y,u,p,q)=-\tfrac{1}{2} u^2 +up +\beta q
\end{align}
The maximum of $u \mapsto H_0 (u)$ is obtained when $\frac{\partial H_0}{\partial u}= -u + p=0$, i.e. when
\begin{equation}
    u=p.
\end{equation}

The adjoint equation is
\begin{align}
    p(dz)&= q(z) B(dz);\quad z < Z=(T,X)\nonumber\\
    p(T,X)&= \theta Y(T,X)
\end{align}
 
Let us try to put
\begin{align}\label{4.18}
    p(t,x)=\lambda(t,x) Y(t,x)
\end{align}
for some deterministic function $\lambda$.
Then by the It\^{o} formula
\begin{align} \label{4.19}
    \lambda(t,x)Y(t,x)&=\lambda(0,0) Y(0,0) \nonumber\\
    &+\int_0^t \int_0^x \Big\{Y(s,a) \frac{\partial^2 \lambda}{\partial s \partial x}(s,a) 
    +\lambda(s,a) u(s,a)\Big\} ds da \\&+ \text{ terms containing } B(ds,da).\nonumber
\end{align}

Using the concept of quadratic variation of $2-$%
parameter martingales (see e.g. Imkeller \cite{I2}), one finds that the decomposition of a (continuous) $2-$%
parameter "semimartingale", which is given by a sum of a $2-$%
parameter process of bounded variation and a $2-$%
parameter martingale, is unique. So, comparing the latter equation with the adjoint equation, we see that we must have 
\begin{align}
 Y(t,x) \frac{\partial^2 \lambda}{\partial t \partial x}(t,x) 
    +\lambda(t,x) u(t,x)=0 \text{ for all } t, x.   
\end{align}
Combining this with \eqref{4.18} we get
\begin{align}
    Y(t,x)\Big[\frac{\partial^2 \lambda}{\partial t \partial x}(t,x) 
    +\lambda^2(t,x)\Big]=0,
\end{align}
with terminal condition
\begin{align}
    \lambda(T,X)=\theta.
\end{align}
 
In addition we get from \eqref{4.18} the other boundary condition
\begin{align}
 \lambda(0,0)=\frac{E[\theta Y(T,X)]}{Y(0,0)}.   
\end{align}
With this choice of $u,p,\lambda$ we see that all the conditions of the sufficient maximum principle are satisfied, and we have proved the following:
\begin{theorem}
 The optimal control $\widehat{u}$ for the LQ problem \eqref{4.3} is given in feedback form by
 \begin{align}
     \widehat{u}(t,x)=\lambda(t,x) Y(t,x); \quad t\leq T,x\leq X,
 \end{align}
 where $\lambda(t,x)$ solves the time-space Riccati equation
 \begin{align}
     \begin{cases}
      \frac{\partial^2 \lambda}{\partial t \partial x}(t,x) 
    +\lambda^2(t,x) = 0; \quad 0 \leq t\leq T, 0\leq x\leq X,\\
    \lambda(T,X)=\theta,\label{Riccati}\\
    \lambda(0,0)=\frac{E[\theta Y(T,X)]}{Y(0,0)}.
     \end{cases}
 \end{align}
\end{theorem}

\begin{remark}
Let $\varphi _{1}$ be a solution to the Riccati equation%
\begin{equation*}
\dot{\varphi}_{1}(t)=(\varphi _{1}(t))^{2},\varphi _{1}(0)=1,0\leq t\leq T
\end{equation*}%
and $\varphi _{2}$ be a solution to%
\begin{equation*}
\dot{\varphi}_{2}(x)=-(\varphi _{2}(x))^{2},\varphi _{2}(0)=\theta ,0\leq
x\leq X.
\end{equation*}%
Define $\lambda (t,x)=$ $\alpha _{1}(t)\alpha _{2}(x)$, where $\alpha
_{1}(t):=\varphi _{1}(T-t)$ and $\alpha _{2}(x):=\varphi _{2}(X-x)$. Then%
\begin{equation*}
\frac{\partial ^{2}\lambda (t,x)}{\partial t\partial x}=\dot{\varphi}%
_{1}(T-t)\dot{\varphi}_{2}(X-x)=(\alpha _{1}(t))^{2}(-(\alpha
_{2}(t))^{2})=-(\lambda (t,x))^{2}
\end{equation*}%
with $\lambda (T,X)=\varphi _{1}(0)\varphi _{2}(0)=\theta $.
\end{remark}

By
solving the Riccati equations, we find that $\lambda $ given by%
\begin{equation*}
\lambda (t,x)=\frac{1}{(1-T+t)(\theta ^{-1}+X-x)}
\end{equation*}%
is an explicit solution to the above hyperbolic PDE with boundary condition $%
\lambda (T,X)=\theta $ for $0<T<1$. 
%\end{remark}
 
Let us now have a look at the other
condition 
$$\lambda (0,0)=\frac{\theta E\left[ Y(T,X)\right]} {Y(0,0)}$$

We observe
that%
\begin{equation*}
E\left[ Y(t,x)\right] =Y(0,0)+\int_{0}^{t}\int_{0}^{x}\lambda (s,a)E\left[
Y(s,a)\right] dsda\text{.}
\end{equation*}%
So, if we use Picard iteration combined with the fact that $\lambda$ can be written as a product of a function in $t$ and another function in $x$, we see that the solution of the latter equation has the
representation%
\begin{equation*}
E\left[ Y(t,x)\right] =Y(0,0)f\left(\int_{0}^{t}\int_{0}^{x}\lambda (s,a)dsda\right),
\end{equation*}%
where the function  $f:\mathbb{R}\mapsto \mathbb{R}$ is defined by 
\begin{equation*}
f(y)=\sum_{n\geq 0}\frac{y^{n}}{(n!)^{2}}\text{.}
\end{equation*}%
 
On the other hand,%
\begin{eqnarray*}
\int_{0}^{T}\int_{0}^{X}\lambda (s,a)dsda &=&\int_{0}^{T}\int_{0}^{X}\frac{1%
}{(1-T+s)(\theta ^{-1}+X-a)}dsda \\
&=&-\log (1-T)\log (1+X\theta ).
\end{eqnarray*}%
So the condition
 $$\lambda (0,0)=\frac{\theta E\left[ Y(T,X)\right]} {Y(0,0)}$$

 is
equivalent to%
\begin{equation*}
\frac{1}{(1-T)(\theta ^{-1}+X)}=\theta f(-\log (1-T)\log (1+X\theta ))
\end{equation*}%
or%
\begin{equation}
1=(1-T)(1+X\theta )f(-\log (1-T)\log (1+X\theta )).  \label{Condition}
\end{equation}%

For given $T<1$ and $\theta >0$ the expression on the right hand side of the
latter equation converges to $(1-T$) for $X\rightarrow 0$. \ For $%
X\rightarrow \infty $, this expression converges to $\infty $. Because of
continuity, we then see that there exists a $X=X(T,\theta )>0$ such that the
equation (\ref{Condition}) is satisfied. Using such a time horizon $X$,
gives the other boundary condition.     
%\end{remark}

\subsection{Example related to machine learning}

In machine learning the (continuous-time) stochastic gradient descent method
(see e.g. \cite{MS} and the references therein) is used to minimize an
objective function $f:\mathbb{R}^{d}\longrightarrow \mathbb{R}$. Compared to the classical gradient descent method without noise, this approach is
especially computationally efficient, when the dimension $d$ in practical optimization problems is high. If the objective
function is sufficiently smooth the critical points of $f$ with respect to
local or global minima may be found by means of solutions to SDEs of the
type
\begin{equation}
dY_{t}=-\eta \nabla f(Y_{t})dt+\beta_0 dB_{t},Y_{0}=x\in \mathbb{R}^{d},t\geq
0\text{,}  \label{MLSDE}
\end{equation}
 
where $\eta \geq 0$ is the learning rate (or step size), $\beta_0 \in \mathbb{%
R}^{d\times d}$, $B_{t},t\geq 0$ a Brownian motion and where $\nabla $
denotes the gradient of a function. In general, the selection of an
"optimal" learning rate $\eta $, which determines the optimal step size
towards a minimum in the sense of speed, is in general difficult. If $\eta $
is chosen too small, the solution may converge too slowly to a critical point.
On the other hand, a too large $\eta $ could result in overshoot or
divergence. In order to gain a deeper understanding of the latter problem,
one may consider instead of the SDE (\ref{MLSDE}) a more general framework
(at the possible expense of computational cost) in connection with the
following type of hyperbolic SPDE:
 
\small
\begin{equation}
Y(t,x)=y-\int_{0}^{t}\int_{0}^{x}u(s,a)\nabla f(Y(s,a))dsda+\beta_0
B(t,x),y\in \mathbb{R}^{d},t,x\geq 0\text{,}  \label{MLHSPDE}
\end{equation}
where $u:\Omega \times \left[ 0,\infty \right) ^{2}\longrightarrow \left[
0,\infty \right) $ is a stochastic learning rate in time and space given by
an adapted random field and where $B$ is a Brownian sheet in $\mathbb{R}^{d}$%
. Formally, by choosing in (\ref{MLHSPDE}) $u=\eta \delta _{x}$ for the
Dirac delta function $\delta _{x}$ in a fixed point $x$ and $\eta \geq
0$ we obtain an SDE of the type (\ref{MLSDE}). So the random field dynamics (\ref{MLHSPDE}) provides a more general framework than that in the one-parameter case (\ref{MLSDE}) for finding the critical points of $f$. 
 
On the other hand, we may
view the integral term%
\begin{equation*}
\int_{0}^{x}u(s,a)\nabla f(Y(s,a))da
\end{equation*}%
in (\ref{MLHSPDE}) for a fixed $x$ and a certain class of  stochastic $2-$%
parameter learning rate processes as an (weighted) average of $\nabla f(Y(s,a))$, $0 \leq a \leq x$ in
(\ref{MLSDE}). Here $Y(s,a), 0 \leq a \leq x$ can be interpreted as a group of mountain hikers in the optimisation landscape who communicate with each other with respect to (average) gradient information in order to find the descent to the valley (i.e. minimum). 
 
The latter, combined with the "exploration ability"
of the Brownian sheet with respect to the spatial parameter direction in the
optimisation landscape, suggests a solution that converges to rather flat
minima, while escaping from sharp minima. The convergence to flat minima, however, is in many applications a
favorable feature from a machine learning point of view (see \cite{HS}). 

In order to construct optimal stochastic $2-$parameter learning rate
processes one may e.g study stochastic control problems based on the
stochastic maximum principle for SPDEs driven by a Brownian sheet with
respect to certain performance functionals as e.g.

\begin{equation}
J(u)=-E\left[ \int_{0}^{T}\int_{0}^{X}u^{2}(s,a)dsda+f(Y(T,X))\right] \text{,} \label{EF}
\end{equation}%
where one mimimises the expected value of $f(Y(T,X))$, while the "energy invested" in $u$ is kept minimal.
 
Using the first order Taylor expansion, we can also approximate $\nabla f$
in (\ref{MLHSPDE}) by an affine function $g$ given by $g(x)=a+Ax$ for $a\in 
\mathbb{R}^{d}$, $A\in 
\mathbb{R}^{d\times d}$ and obtain a more simplified framework for our
stochastic control problem with respect to $u$. In this setting, let us now
consider the case $d=1$ and the following controlled process:             

\begin{align*}
Y_{u}(t,x)  & =Y(0,0)-\int_0^t\int_0^{x}u(s,a)Y_u(s,a)dsda+\int_0^t \int_0^{x}\beta_{0}B(ds,da).\nonumber
\end{align*}
\begin{problem} \label{p3}
    We want to maximise the performance functional
\begin{equation*}
J(u)=-E\Big[\int_{0}^{T}\int_{0}^{X}u^{2}(s,a)dsda+\theta Y^{2}(T ,X)\Big].
\end{equation*}
\end{problem}
 
In this case the associated Hamiltonian is
\begin{align*}
H(t,x,y,u,p,q,L)=-u^2 -yup +\beta_0 q- ( L \star (uy)),
\end{align*}
which gives
\begin{align}
 \frac{\partial H}{\partial y}= -up-(L\star u)   
\end{align}
and hence the adjoint variable $L$ is given by
\begin{align}
    L=-\frac{\partial H}{\partial y}=up+(L\star u) \label{L2}
\end{align}
which is an integral equation for $L(z)$.

The adjoint BSDE is 
 \begin{align}
   p(z)&=p(0)+ \int_{R_{z}} (u(\zeta) p(\zeta)+(L \star u))(\zeta) d\zeta \nonumber\\
   &-\int_{R_{z}} q(\zeta)B(d\zeta)-\iint\limits_{R_{z}\times R_{z}} r(\zeta,\zeta')B(d\zeta)B(d\zeta'),\quad z\leq (T,X)\label{SPDE2a}\\
   p(T,X)&=-2\theta Y(T,X) .
    \end{align}
In order to maximize $H$ with respect to $u$, we first compute
$$
\frac{\partial H}{\partial u}=-2u - yp - ( L \star y)
$$
which gives the optimal control 
$$u=u^{*}=\tfrac{1}{2}(yp+(L\star y))$$
 
We have proved:
\begin{theorem}
The optimal control $u(z)=u^{*}(x)$ for Problem \ref{p3} is 
given as the 4th component of the solution $(Y(z),L(z),p(z),u(z))$ of the coupled system of equations, consisting of
\begin{align}
    u(z)&=u^{*}(z)=\tfrac{1}{2}[Y(z)p(z)+(L\star Y)(z)],\\
Y(t,x)& =Y(0,0)-\int_0^t\int_0^{x}u(s,a)Y_u(s,a)dsda+\int_0^t \int_0^{x}\beta_{0}B(ds,da),\\
L(z)&= u(z)p(z) +(L\star u)(z)
\end{align}
\end{theorem}

and 
\begin{align}
\begin{cases}
p(z)&=p(0)+ \int_{R_{z}} (u(\zeta) p(\zeta)+(L \star u))(\zeta) d\zeta \\
   &-\int_{R_{z}} q(\zeta)B(d\zeta)-\iint\limits_{R_{z}\times R_{z}} r(\zeta,\zeta')B(d\zeta)B(d\zeta'),\quad z\leq (T,X)\label{spde2}\\
   p(T,X)&=-2\theta Y(T,X) .
   \end{cases}
    \end{align}
%\end{theorem}
 
\begin{remark}
\bigskip In the more general case, when $\triangledown f(x)=Ax$ for $A\in 
\mathbb{R}^{d\times d}$, one can show a similar result about the optimal control $u^{\ast }$ with
respect to the controlled process (\ref{MLHSPDE}) and performance functional
(\ref{EF}) 
%is given by $u^{\ast }(t,x)=-\frac{1}{2}(\triangledown
%f(Y(t,x)))^{\ast }p(t,x)$, where $(Y,p)$ solves a corresponding
%forward-backward system of SPDEs ($^{\ast }$ denotes transpose).
\end{remark}
%%%%%%%%%%%%%%%%%%%%%%%%%%%%%%%%%%%%%%%%%%%%%%%%%%%%%%%%%%%%%%%%%%%%%%%%%%%%%%%%%%%%%%%%%%%%%%%%%%%%%%%%%%%%%%%%%%%%%%%%%%%%%%%%%%%%%%%%%%%%%%%%%%%%%%

\end{document}